\documentclass[11pt]{article}

\usepackage{amsthm,amssymb}
\catcode`\@=11
\@addtoreset{equation}{section}

\newtheorem{thm}{Theorem}[section]
\newtheorem{prop}{Proposition}[section]
\newtheorem{cor}{Corollary}[section]
\newtheorem{rem}{Remark}[section]
\newtheorem{lem}{Lemma}[section]
\newtheorem{dfn}{Definition}[section]

\newcommand{\eps}{\varepsilon}

\newcommand{\sgn}{{\mathrm{sgn}\,}}

\newcommand{\codim}{{\mathrm{codim}\,}}

\newcommand{\id}{{\mathrm{id}}}

\newcommand{\R}{{\mathbb{R}}}
\newcommand{\bbf}{{\mathbf{b}}}
\newcommand{\cbf}{{\mathbf{c}}}
\newcommand{\kbf}{{\mathbf{k}}}
\newcommand{\hbf}{{\mathbf{h}}}

\newcommand{\nbf}{{\mathbf{n}}}
\newcommand{\xbf}{{\mathbf{x}}}
\newcommand{\ybf}{{\mathbf{y}}}
\newcommand{\zbf}{{\mathbf{z}}}

\newcommand{\tbf}{{\mathbf{t}}}
\newcommand{\sbf}{{\mathbf{s}}}
\newcommand{\qbf}{{\mathbf{q}}}
\newcommand{\ubf}{{\mathbf{u}}}
\newcommand{\vbf}{{\mathbf{v}}}
\newcommand{\wbf}{{\mathbf{w}}}

\newcommand{\Hbf}{{\mathbf{H}}}
\newcommand{\C}{{\mathbb{C}}}
\newcommand{\Z}{{\mathbb{Z}}}
\newcommand{\T}{{\mathbb{T}}}

\newcommand{\calA}{\mathcal{A}}

\newcommand{\calH}{\mathcal{H}}
\newcommand{\calL}{\mathcal{L}}
\newcommand{\calM}{\mathcal{M}}
\newcommand{\calN}{\mathcal{N}}
\newcommand{\calF}{\mathcal{F}}
\newcommand{\calP}{\mathcal{P}}
\newcommand{\calQ}{\mathcal{Q}}
\newcommand{\calR}{\mathcal{R}}
\newcommand{\calS}{\mathcal{S}}
\newcommand{\calT}{\mathcal{T}}

\newcommand{\Real}{{\mathrm{Re}\,}}

\title{Degenerate billiards
\author{Sergey Bolotin\\
Moscow Steklov Mathematical Institute\\
and\\
University of Wisconsin-Madison
} }

\begin{document}

\maketitle

\begin{abstract}
In an ordinary billiard system trajectories of a Hamiltonian system are
elastically reflected after a collision with a hypersurface (scatterer).
If the scatterer is a submanifold of codimension more than one,
we say that the billiard is degenerate. Then collisions are rare.
We study trajectories   of  degenerate billiards which have an infinite number of
collisions with the scatterer.
Degenerate billiards appear as limits of systems with elastic reflections or as limits of systems with singularities in celestial mechanics. We prove the existence of trajectories  of such systems
shadowing trajectories of the corresponding degenerate billiards.
The proofs are based on a version of the method of an anti-integrable limit.
\end{abstract}

\section{Degenerate billiards}

\subsection{Definitions}

\label{sec:bill}

Consider  a Hamiltonian system $(M,H)$ with the configuration space $M$ and smooth (at least $C^3$)
Hamiltonian $H(q,p)$. The symplectic structure $\omega=dp\wedge dq$  on the phase space  $T^*M$ is standard, so we do not include it in the notation.\footnote{Locally a twisted symplectic structure  on $T^*M$ can be replaced
 by a standard one by a calibration transformation. All results in this paper are essentially local.} Suppose that $H$ is
strictly convex and superlinear in the momentum $p\in T_q^*M$. Then a solution $(q(t),p(t))\in T^*M$ is
determined by the trajectory $q(t)\in M$ of the corresponding Lagrangian system
with the Lagrangian
$$
L(q,\dot q)=\max_p (\langle p,\dot q\rangle-H(q,p)).
$$

Let $N\subset M$ be a submanifold in $M$ which is called a scatterer.  Suppose that when a trajectory
$q(t)$ of  system $(M,H)$ meets the scatterer at a collision point
$x=q(\tau)\in N$, it is reflected according to the elastic reflection law\footnote{Here $p_-$ is the momentum after the collision,
and $p_+$ before the collision. The strange notation is chosen to fit with the notation for the initial and final momenta of a collision orbit $\gamma:[t_-,t_+]\to M$.}
\begin{eqnarray}
\label{eq:Deltap} \Delta p(\tau)=p_--p_+\perp T_x N,\qquad p_\pm=p(\tau\mp 0),\\
\Delta H(\tau)=H(x,p_-)-H(x,p_+)=0.
\label{eq:DeltaH}
\end{eqnarray}
Thus the tangent component $y\in T_x^*N$ of the momentum $p\in T_x^*M$ is preserved.
We always assume that the momentum has a jump at the collision: $\Delta p(\tau)\ne 0$.
If $p_+=p_-$,
the trajectory goes through $N$ without noticing.

When $N$ is a hypersurface bounding a domain $D$ in $M$,   we obtain
a usual billiard system $(D,N,H)$ in $D$ with boundary $N$.  In this paper we consider the case when the codimension of $N$ in $M$ is greater than 1. Then we say that $(M,N,H)$ is a  degenerate billiard.

Trajectories of the  degenerate billiard having collisions with
$N$ form a zero measure set  in the phase space.
Moreover  $p_+$ does not determine $p_-$ uniquely by
(\ref{eq:Deltap})--(\ref{eq:DeltaH}): for given $p_+$ the set of possible $p_-$
has dimension $\codim N-1$.
A degenerate billiard  is not a dynamical system according to the usual definition:
the past of a trajectory does not determine its future after a collision.
The simplest example is a 0-dimensional billiard with $N$ a discrete set in $M$.
Then only  continuity of the Hamiltonian (\ref{eq:DeltaH}) remains, so a trajectory can be reflected in any direction
after a collision.

We are interested in trajectories   with multiple collisions
which are called collision chains.  Collision
chains $\gamma:[\alpha,\beta]\to M$  are extremals of the
action functional
\begin{equation}
\label{eq:I}
I(\gamma)=\int_\alpha^\beta L(\gamma(t),\dot\gamma(t))\,dt=\sum_{j=0}^nI(\gamma_j),\qquad \gamma_j=\gamma|_{[t_j,t_{j+1}]}
\end{equation}
on the set of curves $\gamma:[\alpha,\beta]\to M$ with fixed end points $a=\gamma(\alpha)$ and $b=\gamma(\beta)$, subject to the
constraints $\gamma(t_j)= x_j\in N$ for some sequence $\alpha=t_0<t_1<\dots <t_n<t_{n+1}=\beta$. Collision points
$x_j\in N$ and collision times $t_j\in (\alpha,\beta)$   are independent variables.

Indeed, if
$\gamma$ is an extremal of $I$, then each segment
$\gamma_j=\gamma|_{[t_j,t_{j+1}]}$ is a trajectory of the
Hamiltonian system and by the first variation formula
$$
\delta I(\gamma)=\sum_{j=1}^n (\Delta p(t_j)\,\delta x_j-\Delta H(t_j)\,\delta t_j)=0,\qquad \delta
x_j\in T_{x_j}N,
$$
which implies
\begin{equation}
\label{eq:elastic2}
\Delta p(t_j)\perp T_{x_j} N,\quad\Delta H(t_j)=0.
\end{equation}

We call a  trajectory  $\gamma:[t_-,t_+]\to M$ of system $(M,H)$ a collision orbit if
its end points lie in  $N$ and there is no tangency and no early collisions with the scatterer:
\begin{equation}
\label{eq:nocoll}
\gamma(t_\pm)=x_\pm\in N,\qquad
v_\pm=\dot\gamma(t_\pm)\notin T_{x_\pm}N,\qquad
\gamma(t)\notin N\quad\forall\; t\in (t_-,t_+).
\end{equation}
An infinite collision chain $\gamma:\R\to M$
is a concatenation  of a sequence $\gamma=(\gamma_j)_{j\in\Z}$ of
collision orbits $\gamma_j:[t_j,t_{j+1}]\to M$  joining the points $x_j,x_{j+1}\in N$ such that the
elastic reflection law  (\ref{eq:elastic2})  is satisfied at each
collision.
Equivalently, $\gamma$ is an extremal of the formal action functional
$$
I(\gamma)=\int_{-\infty}^\infty L(\gamma(t),\dot\gamma(t))\,dt=\sum_{j\in\Z} I(\gamma_j)
$$
on the set of curves with the constraints
$\gamma(t_j)= x_j\in N$
for variations with finite support. We also require the  jump condition  $\Delta p(t_j)\ne 0$.

The Hamiltonian is constant along a collision chain, so let us fix energy.
The  restriction of the Hamiltonian system to the energy level $H=E$ will be denoted $(M,H=E)$.
 Trajectories $\gamma:[\alpha,\beta]\to M$ of system  $(M,H=E)$  are extremals of the Maupertuis action
$$
J_E(\gamma)=\int_\alpha^\beta \|\dot\gamma(t)\|_E\,dt,
$$
i.e.\ geodesics of the Jacobi metric
\begin{equation}
\label{eq:metric}
\|\dot q\|_E=\max_p \{\langle p,\dot q\rangle:H(q,p)=E\}
\end{equation}
in the domain of possible motion
\begin{equation}
\label{eq:DE}
D_E=\{q\in M:W(q)<E\},\qquad W(q)=\min_pH(q,p)=-L(q,0).
\end{equation}
For trajectories with energy $E$,
$$
J_E(\gamma)=\int_\gamma p\, dq.
$$

\begin{rem}
In general the metric is convex but not positive definite in $D_E$, so it is not a Finsler metric.
However, in a neighborhood of any point  $q_0\in D_E$ it can be assumed to be positive definite. Indeed,
 suppose for $q=q_0$ the maximum in (\ref{eq:metric}) is attained at $p=p_0$.
Then  the calibration transformation $p\to p-\nabla f(q)$, $\nabla f(q_0)=p_0$, does not change the trajectories in $M$
but now $H(q_0,0)=W(q_0)$. Then the metric is positive definite near $q_0$.
\end{rem}

In applications, $H$ is usually quadratic in the momentum:
\begin{equation}
\label{eq:HL}
H(q,p)=\frac12\|p-w(q)\|^2+ W(q),\quad L(q,\dot q)=\frac12\|\dot q\|^2+\langle w(q),\dot q\rangle -W(q),
\end{equation}
where $\|\;\|$ is a Riemannian metric on $M$, and $w$ -- a covector field.\footnote{We use the same notation $\|\; \|$ for the norm of a vector and a covector.}
For the classical Hamiltonian (\ref{eq:HL}),
$$
E=\frac12\|\dot q\|^2+W(q),\qquad \|\dot q\|_E=\sqrt{2(E-W(q))}\|\dot q\|+\langle w(q),\dot q\rangle.
$$

When the energy $E$ is fixed, we denote the degenerate  billiard by $(M,N,H=E)$.
Its collision chains  are extremals of the formal Maupertuis action functional
$$
J_E(\gamma)=\sum_{j\in \Z} J_E(\gamma_j)
$$
on the set of chains $\gamma=(\gamma_j)_{j\in\Z}$ of nonparametrized curves  joining a sequence of points  $x_j\in N$.

\subsection{Collision map}

\label{sec:coll}

As was already mentioned, a degenerate billiard is not a dynamical system:
the past of a trajectory does not determine its future after the collision.
Let us define a discrete dynamical system  whose trajectories correspond to  collision chains.

Fix energy $H=E$.
 Let $\gamma:[t_-,t_+]\to D_E$ be a collision orbit joining the points $a_-,a_+\in N$.
We call $\gamma$ nondegenerate if  $a_-$ and $a_+$ are non-conjugate along the extremal $\gamma$ of the functional $J_E$.
Then there exists a neighborhood  $U\subset M\times M$ of $(a_-,a_+)$ such that for all
$(q_-,q_+)\in U$ there exists an   orbit $\gamma(q_-,q_+)$ with energy $E$ joining $q_-$ and $q_+$, and it smoothly depends on $q_-,q_+$.
 The action
$$
S(q_-,q_+)=J_E(\gamma(q_-,q_+))
$$
is a smooth function on $U$. The initial and final momenta of the orbit $\gamma$  are
$$
p_-=-D_{q_-}S,\quad p_+=D_{ q_+}S.
$$
The action function  satisfies the Hamilton-Jacobi equation
$$
H(q_-, -D_{q_-}S)=H(q_+, D_{ q_+}S)=E.
$$
The twist of the action function  is the linear transformation
$$
B(q_-,q_+)= D_{ q_-}D_{ q_+}S:T_{q_-}M\to T_{q_+}^*M,
$$
i.e.\ a bilinear form    on $T_{q_-}M\times T_{q_+}M$.
Differentiating the Hamilton-Jacobi equation, we obtain
\begin{equation}
\label{eq:twist}
B(a_-,a_+)v_-=0,\quad B^*(a_-,a_+)v_+=0,\qquad v_\pm=\dot \gamma(t_\pm).
\end{equation}
Thus for autonomous systems,  $B(a_-,a_+)$ is always degenerate.

 We say that the collision orbit $\gamma$ has nondegenerate
twist if the restriction
of the bilinear form $B(a_-,a_+)$ to $T_{a_-}N\times T_{a_-}N$    is nondegenerate.
For  this it is necessary that $v_\pm\notin
T_{a_\pm}N$, i.e.\ the collision orbit is not tangent to the scatterer $N$ at the end points.
For an ordinary billiard, when $N$ is a hypersurface, this is also sufficient for the nondegenerate twist,
but not for a degenerate billiard.

If $\gamma$ has nondegenerate  twist, the restriction of $S$   to a neighborhood
of $(a_-,a_+)$ in $N\times N$ is the generating function of a locally
defined symplectic map $f:V^-\to V^+$ of open sets
$V^\pm\subset T^*N$:
$$
f(x_{-},y_{-})=(x_{+},y_{+})\quad \Leftrightarrow\quad
y_+= D_{ x_+}S,\quad y_-=-D_{ x_-}S.
$$
Here $y_{\pm}\in T_{x_\pm}^*N$ are the   tangent  projections of the
collision momenta $p_\pm\in T_{x_\pm}^*M$ of the collision orbit $\gamma(x_-,x_+)$.

In general there may exist several (or none) nondegenerate
collision orbits  with   energy $E$ joining a
pair of points in $N$. Thus we obtain a collection   $\calL=\{ L_k\}_{k\in K}$ of action functions on
open sets $U_k\subset N\times N$. Under the twist condition,    $L_k$ generates a
local symplectic map $f_k:V_k^-\to V_k^+$ of open sets in $T^*N$.  Let $\calF=\{f_k\}_{k\in K}$.
We call the partly defined multivalued map\footnote{It is more correct to call $\calF$  a relation.} $\calF:T^*N\to T^*N$ the collision map,
or the scattering map of the degenerate billiard.
It defines a discrete dynamical system    -- the skew product of the maps
$\calF=\{f_k\}_{k\in K}$ which is a map of a subset in $T^*N\times K^\Z$.

An orbit of $\calF$ is a pair $(\kbf,\zbf)$ of sequences $\kbf=(k_j)$,
$\zbf=(z_j)$, where $z_j=(x_j,y_j) \in
V_{k_j}^-\cap V_{k_{j-1}}^+$, such that $z_{j+1}=f_{k_j}(z_j)$.  The orbit
$(\kbf,\zbf)$ defines  a chain  of collision orbits
$\gamma_j$ joining $x_j$ with $x_{j+1}$.  The tangent collision  momenta of the collision chain are
\begin{equation}
\label{eq:y_j}
y_j=  D_{x_j} L_{k_{j-1}}(x_{j-1},x_{j})=-D_{x_j} L_{k_j}(x_j,x_{j+1})\in T_{x_j}^*N.
\end{equation}

Not all orbits of $\calF$  define admissible collision chains: we need also the  jump condition $\Delta p_j\ne 0$.

\subsection{Zero-dimensional  billiards}

\label{sec:zero}

Suppose the scatterer $N=\{a_i\}_{i\in J}$ is a finite set in $M$. Fix  energy $E>\max_NW$ and let $\{\gamma_k\}_{k\in K}$ be a finite collection of nondegenerate  collision orbits of the degenerate billiard $(M,N,H=E)$.
Denote by $a_k^\pm\in N$ the initial and final points of $\gamma_k$
  and by $p_k^\pm$ the initial and final momentum.  Consider an oriented graph
  with the set of vertices $K$  and the set of edges $\Gamma$. Two vertices $k,k'\in K$ are joined by an edge $(k,k')\in \Gamma$ if $a_k^+=a_{k'}^-$ and $p_k^+\ne p_{k'}^-$ (jump condition).
Then to any path $\kbf=(k_j)_{j\in\Z}$ in the   graph $\Gamma$ there corresponds a collision  chain
 $\gamma=(\gamma_{k_j})_{j\in\Z}$ of the degenerate billiard $(M,N,H=E)$.

 Let $\Sigma_\Gamma\subset K^\Z$ be the set of all  paths $\kbf$  in the graph.
 Then dynamics of the degenerate billiard $(M,N,H=E)$ is conjugate to a  topological Markov chain
 (called also a subshift of finite type \cite{K-H}) $\calT:\Sigma_\Gamma\to\Sigma_\Gamma$,
 $(k_j)\to (k_{j+1})$.
 The collision map $\calF$ is a multivalued map of a finite set $K$.
 Such  trivial billiards appear in nontrivial applications, see sections \ref{sec:small} and \ref{sec:cent}.

\subsection{Discrete action functional}

\label{sec:func}

In this paper we avoid using the collision map. Indeed,
in general the twist condition is not satisfied and moreover it is not easy to check.
In some application the scatterer $N$ has connected components  with different dimensions.
Then there is no chance for the twist condition to hold.

Without the twist condition,
the degenerate billiard $(M,N,H=E)$ can be viewed as a discrete Lagrangian
system (DLS) with multivalued Lagrangian $\calL=\{ L_k\}_{k\in K}$, see \cite{Bol-Tre:ANTI}.
Collision chains with energy $E$  correspond to  critical points
$\xbf=(x_j)_{j\in\Z}$ of the formal discrete action functional
\begin{equation}
\label{eq:Ak}
A_\kbf(\xbf)=\sum_{j\in\Z}  L_{k_j}(x_j,x_{j+1}),\qquad (x_j,x_{j+1})\in U_{k_j},
\end{equation}
where $k_j\in K$ can be chosen randomly at each step. Thus a trajectory of the DLS is a pair $(\kbf,\xbf)\in K^\Z\times N^\Z$
such that $A'_\kbf(\xbf)=0$.

For infinite collision chains, the   sum  (\ref{eq:Ak}) makes no sense, so $A_\kbf(\xbf)$ is a
formal functional,  but $D_{x_j} A_\kbf(\xbf)\in T_{x_j}^*N$ is well defined.
Thus the derivative and the Hessian
$$
A'_\kbf(\xbf)=(D_{x_j} A_\kbf(\xbf))_{j\in\Z},\qquad A''_\kbf(\xbf)=(D_{x_i}D_{x_j} A_\kbf(\xbf))_{i,j\in\Z}
$$
are well defined. If we fix a Riemannian metric and identify $T_{x_i}N$ and $T_{x_i}^*N$,
the Hessian becomes a linear operator
$A''_\kbf(\xbf):l_\infty \to
l_\infty $, where $l_\infty$ is the set of sequences
$$
\ubf=(u_i)_{i\in \Z},\qquad u_i\in T_{x_i}N,\quad \|\ubf\|_\infty=\sup_i\|u_i\|<\infty.
$$

The Hessian operator  is 3-diagonal:
\begin{equation}
\label{eq:var}
A_\kbf''(\xbf)\ubf= \vbf,\qquad v_i=B_{i-1}u_{i-1} + A_i u_i + B_i^* u_{i+1},
\end{equation}
where $A_i,B_i$   are linear operators
$$
A_i=D^2_{x_i}(L_{k_{i-1}}(x_{i-1},x_{i})+L_{k_i}(x_i,x_{i+1})),\qquad  B_i=D_{x_i}D_{x_{i+1}}L_{k_i}(x_i,x_{i+1}).
$$
If $K$ is finite and for each $k\in K$ there is a compact set $X_k\subset U_k$ such that $(x_i,x_{i+1})\in X_{k_i}$ for all $i$, then $A_\kbf''(\xbf):l_\infty\to l_\infty$ is a bounded operator.

The variational equation of the trajectory $(\kbf,\xbf)$ is $A_\kbf''(\xbf)\ubf=0$.
Under the twist condition, $B_i$ is invertible, and the variational equation defines the linear
Poincar\'e map $P_i:(u_{i-1},u_i)\to (u_i,u_{i+1})$.

Without the twist condition, the dynamics of the DLS
is represented by a translation $\calT:K^\Z\times N^\Z\to K^\Z\times N^\Z$, $(k_j,x_j)\to (k_{j+1},x_{j+1})$.
We equip $K^\Z\times N^\Z$ with a product topology.
If $\calT$ has  a compact invariant set $\Lambda\subset K^\Z\times N^\Z$ of  trajectories of the DLS, and the
collision map $\calF$ is well defined, then it will have a compact invariant set $\Lambda_\calF\subset K^\Z\times N^\Z$
with $\calF:\Lambda_\calF\to\Lambda_\calF$  topologically conjugate to $\calT:\Lambda\to\Lambda$.

Finite collision chains of the degenerate billiard joining the points $a,b\in M$ are critical points of a finite sum
\begin{equation}
\label{eq:finite}
A_\kbf^{a,b}(\xbf)=\sum_{j=0}^{n}  L_{k_j}(x_j,x_{j+1}),\qquad x_0=a,\quad x_{n+1}=b,\quad \xbf=(x_1,\dots,x_n).
\end{equation}
We call the finite collision chain nondegenerate if the critical point is nondegenerate.

For $n$-periodic collision chains, $\xbf$ is a critical point the periodic action functional
\begin{equation}
\label{eq:per} A_\kbf^{(n)}(\xbf)=\sum_{j=0}^{n-1}
L_{k_j}(x_j,x_{j+1}),\qquad \xbf=(x_1,\dots ,x_n),\quad
x_{n}=x_0 .
\end{equation}
We call the periodic orbit $(\kbf,\xbf)$ nondegenerate if $\xbf$ is a nondegenerate  critical point of  $A_\kbf^{(n)}$.
If the twist condition holds, this is equivalent to the usual nondegeneracy condition  $\det (P-I)\ne 0$, where $P=P_n\circ\dots\circ P_1$ is  the linear Poincar\'e map.

\subsection{Nonautonomous billiards}

\label{sec:time}

If the Hamiltonian $H(q,p,t)$ and the scatterer $N_t$  are time dependent, only some notations need to be changed.
Suppose $N_t=\phi_t(N)$, where $\phi_t:N\to M$ is an embedding. Then the Hamiltonian changes at the collision at $q=\phi_t(x)$ as follows:
$$
\Delta H=\langle\Delta p, D_t\phi_t(x)\rangle,\qquad \Delta p\perp T_q N_t.
$$
Hamilton's variational principle (\ref{eq:I}) is the same.
The discrete Lagrangian  is defined as the action $ L_k(x_-,t_-,x_+,t_+)=I(\gamma)$ of a nondegenerate collision orbit $\gamma:[t_-,t_+]\to M$ with non-conjugate
end points $q_-=\phi_{t_-}(x_-)$ and $q_+=\phi_{t_+}(x_+)$.   Thus $L_k$ is
a function on an open  set $U_k\subset N^2\times \R^2$. Initial and final tangent momenta $y_\pm\in T_{x_\pm}^*N$
and the Hamiltonian  $h_\pm=H|_{t_\pm}$ are given by
\[
y_+= D_{ x_+} L_k,\quad y_-=-D_{ x_-} L_k,\quad h_+= -D_{ t_+} L_k,\quad h_-=D_{ t_-} L_k.
\]
If the twist condition holds, the corresponding collision map will be a symplectic map
of  a subset in $T^*N\times\R^2\{t,h\}$.

Infinite collision chains   correspond to  critical points
$(\xbf,\tbf)=(x_j,t_j)_{j\in\Z}$ of the formal discrete action functional
\begin{equation}
A_\kbf(\xbf,\tbf)=\sum_{j\in\Z}  L_{k_j}(x_j,t_j,x_{j+1},t_{j+1}),\qquad (x_j,t_j,x_{j+1},t_{j+1})\in U_{k_j}.
\end{equation}
A trajectory of the DLS with the Lagrangian
$\calL=\{ L_k\}_{k\in K}$ is a triple $(\kbf,\xbf,\tbf)\in K^\Z\times N^\Z\times\R^\Z$.

Formally we can regard nonautonomous billiard as an autonomous one with one more degree of freedom.
But then the convexity assumptions will not hold.

 Note that already 0-dimensional nonautonomous billiards are  nontrivial, then the collision map is a twist map of $\R^2\{t,h\}$.
Such degenerate billiard appears in the elliptic restricted 3 body problem \cite{Bol:DCDS}.

\subsection{Hyperbolic invariant sets of degenerate billiards}

\label{sec:hyp}

The usual definition of hyperbolicity is formulated in terms of the dichotomy of stable and unstable trajectories
of the variational equation. It works under the twist condition: the operators $B_i$
in (\ref{eq:var}) are invertible, so the linear Poincar\'e maps $P_i$ are well defined.

\begin{dfn}\label{def1}
We say that the trajectory $(\kbf,\xbf)$ of the DLS is hyperbolic if for any $j$ there are stable and unstable subspaces $E_j^\pm(\kbf,\xbf)\subset T_{x_{j}}N\times T_{x_{j+1}}N$  such that
$E_j^+\cap E_j^\pm=\{0\}$ and
$w_j=(u_{j},u_{j+1})\in E^+_j$
 implies $w_i=(u_{i},u_{i+1})\in E^+_i$ for all $i>j$. Moreover $u_i$ decreases exponentially as $i\to\infty$: there is $C>0$ and $\mu\in (0,1)$ such that
 $$
 \|w_i\|\le C\mu^{i-j}\|w_j\|,\qquad i>j.
 $$
 Similarly for the unstable subspace:  $w_j=(u_{j},u_{j+1})\in E^-_j$
 implies $w_i=(u_{i},u_{i+1})\in E^-_i$ for all $i<j$ and $u_i$ decreases exponentially as  $i\to-\infty$:
 $$
 \|w_i\|\le C\mu^{j-i}\|w_j\|,\qquad i<j.
 $$
 We say that a compact $\calT$-invariant set $\Lambda$ of trajectories is hyperbolic if this holds for
 every trajectory $(\kbf,\xbf)\in\Lambda$ with $C,\mu$ independent of the trajectory.
\end{dfn}

 For our purposes another definition, not requiring the twist condition,  is more convenient.

\begin{dfn}\label{def2}
We say that the trajectory $(\kbf,\xbf)$ is hyperbolic if
 the Hessian  $A''_\kbf(\xbf)$  has bounded inverse in the $l_\infty$ norm.
 We say that a compact $\calT$-invariant  set $\Lambda\subset K^\Z\times N^\Z$ of trajectories of the DLS   is hyperbolic
 if this is true for all trajectories:
 $\|(A''_\kbf(\xbf))^{-1}\|_\infty\le C$ with $C$ independent of $(\kbf,\xbf)\in\Lambda$.
 \end{dfn}

If the twist condition holds, so that the collision map is well defined, then, as shown in \cite{Aubry-MacKay}, this definition of hyperbolicity (phonon gap condition) is equivalent to the standard one.

We need a slight modification of this result:

 \begin{prop}
 If $\Lambda$ is compact hyperbolic invariant set, then there is a constant $C>0$ such that
 $$
 \|(A''_\kbf(\xbf))^{-1}\|_\infty\le C
 $$
 for all $(\kbf,\xbf)\in\Lambda$.
\end{prop}

\proof
   Let us show that the operator $\Hbf=A''_\kbf(\xbf)$ is invertible.
We need to solve the equation $\Hbf \vbf=\wbf$.
First suppose that $w_i=0$ for $i\ne j$.
Since $\|v_i\|$ is bounded as $|i|\to\infty$,
necessarily $(v_{j},v_{j+1})\in E_{j}^+$ and $(v_{j-1},v_j)\in E_{j-1}^-$.
Then equation $\Hbf \vbf=\wbf$ reads
$$
B_{j-1}v_{j-1}+B_j^*v_{j+1}+A_jv_j=w_j,\quad
(v_j,v_{j+1})\in E_{j}^+,\quad (v_{j-1},v_j)\in E_{j-1}^-.
$$

For $w_j=0$ this implies $(v_j,v_{j+1})\in E_j^-$. Since $E_j^+\cap E_j^-=\{0\}$, for $w_j=0$ the only solution
of these equations is $v_j=v_{j-1}=v_{j+1}=0$.
Hence for any $w_j$ there exists a unique solution $v_i=G_{ij}w_j$
with the Green function $G_{ij}$ satisfying
$\|G_{ij}\|\le Ce^{-\lambda |i-j|}$,
because trajectories in $E_i^\pm$ tend to 0 exponentially as $|i|\to\infty$.

Now for any $\wbf=(w_i)_{i\in\Z}$, we get formally
$$
v_i=\sum_{j\in\Z}G_{ij}w_j.
$$
If $\|\wbf\|_\infty=\sup_{i\in\Z}\|w_i\|<\infty$, then
$$
\|\vbf\|_\infty=\sup_{i\in\Z}\|v_i\|\le c\|\wbf\|_\infty,\qquad c= C \sum_{i\in\Z}e^{-|i|\lambda}<\infty.
$$
Thus $\Hbf$ is $l_\infty$-invertible.
\qed

\medskip

Under the twist condition, the converse is also true but we do not need this since our work definition will be in terms of the Hessian.

\subsection{Routh reduction of symmetry}

\label{sec:reduce}

If the degenerate billiard has a  first integral, then in general it has no hyperbolic invariant sets or nondegenerate periodic orbits until we reduce symmetry.
We describe the  simplest situation arising in applications. Suppose that there is a one-parameter symmetry group $\Phi_\theta:M\to M$, $\theta\in\R$ or $\theta\in\T=\R/\Z$,
which preserves the Hamiltonian and the scatterer:
$$
\Phi_\theta (N)=N,\qquad H(\Phi_\theta (q),p)=H(q,(D\Phi_\theta(q))^*p).
$$
Let 
$$
u(q)=\frac{\partial}{\partial \theta}\Big|_{\theta=0}\Phi_\theta(q)
$$ 
be the vector field generating $\Phi_\theta$ and $G(q,p)=\langle u(q),p\rangle$ the corresponding Noether integral
of the Hamiltonian system $(M,H)$.
Since $u$ is tangent to $N$, $G$ is preserved by the reflection and so it will be an integral of the degenerate billiard $(M,N,H)$.
The corresponding discrete Lagrangian system $\calL=\{ L_k\}$ will have the symmetry
$$
 L_k(\Phi_\theta (x_-),\Phi_\theta (x_+))= L_k(x_-,x_+).
$$
We can assume that the domain $U_k$ is invariant.
The action functional is also invariant:
$$
A_\kbf(\Phi_\theta \xbf)=A_\kbf(\xbf).
$$
Hence $\ubf=(u(x_j))_{j\in\Z}$ is always in the kernel of the Hessian $A''_\kbf(\xbf)$,
and the Hessian is non-invertible: there are no hyperbolic trajectories except fixed points of $\Phi_\theta$.

For any trajectory $(\kbf,\xbf)$ of the DLS let $y_j$ be   the momentum (\ref{eq:y_j}).
The Noether integral
$$
G_j=\langle u(x_j),y_j\rangle=G
$$
is constant along a trajectory of the DLS.

For fixed $G$, we can reduce symmetry by using the Routh reduction \cite{AKN},
replacing $M$ by $\tilde M=M/\Phi_\theta$
and $N$ by $\tilde N=N/\Phi_\theta$.
This can be also done for the DLS.

We need to assume that the fibration $\pi:N\to \tilde N$ to the orbits of the group action is trivial. This is always true if $\theta\in\R$, and locally if $\theta\in \T$. Then $\tilde N$ can be represented as a cross section $\tilde N\subset N$ of the group action.
For fixed value of the integral $G$, define the reduced Lagrangian (Routh function) by the Legendre transform
\begin{equation}
\label{eq:Routh}
\tilde  L_k(x_-,x_+)= (L_k(x_-,\Phi_\theta(x_+))-G\theta)\big|_{\theta=\theta_k(x_-,x_+)},\qquad x_\pm\in \tilde N,
\end{equation}
where $\theta=\theta_k(x_-,x_+)$ is a  critical  point with respect to $\theta$. This requires nondegeneracy condition:
$$
\langle B_k(x_-,x_+)u(x_-),u(x_+)\rangle\ne 0,
$$
where $B_k$ is the twist of $L_k$. The reduced Lagrangian is locally defined on an open set
$\tilde U_k\subset \tilde N\times \tilde N$.

Then for any trajectory $(\kbf,\xbf)$  with Noether integral $G$ setting  $\tilde \xbf=(\tilde x_j)$, $\tilde x_j=\pi(x_j)$,
we obtain a trajectory $(\kbf,\tilde\xbf)$ of the reduced DLS with the Lagrangian $\tilde\calL=\{\tilde L_k\}_{k\in K}$.

Conversely, a trajectory $(\kbf,\tilde\xbf)$ of the reduced DLS defines a (nonunique) trajectory $(\kbf,\xbf)$ of the original DLS
with Noether integral $G$. Here $x_j=\Phi_{\theta_j}(\tilde x_j)$ for some choice of $\theta_j$.

We need such reduction i.e.\ for billiards appearing in Celestial Mechanics,
see section \ref{sec:cel}.

\section{Degenerate billiards as limits of ordinary billiards}

\subsection{Billiards with small scatterers}

\label{sec:small}

 Let us illustrate  how   0-dimensional billiards appear in real billiard systems. Let $D$ be  a domain in $M$ with smooth boundary, $N=\{a_i\}_{i\in J}$ a finite set in $D$,
  and   $B_i\subset D$  a small convex neighborhood of $a_i$
 with boundary $\partial B_i$. More precisely, let $f_i:T_{a_i}M\to M$ be a smooth map such that $f_i(0)=a_i$ and $Df_i(0)=\id$
 (for example, the exponential map $f_i(u)=\exp_{a_i}u$).  We assume that $\partial B_i=f_i(\eps S_i)$,
 where $S_i$ is a strictly convex hypersurface (with positive definite second fundamental form) in $T_{a_i}M$ containing 0. Let
 $$
 \Omega_\eps=D\setminus \bigcup_{i\in J} B_i,\qquad
 \Sigma_\eps=\partial \Omega_\eps=\bigcup_{i\in J}\partial B_i\cup\partial D.
 $$
  Consider a usual billiard system $(\Omega_\eps,\Sigma_\eps,H=E)$  in the domain $\Omega_\eps$ with the boundary $\Sigma_\eps$.
It is  natural to expect that as $\eps\to 0$ its trajectories will approach trajectories of the   degenerate billiard $(D,N\cup \partial D,H=E)$. Note that the scatterer consists of several components, 0-dimensional set $N$ and the hypersurface $\partial D$.

More precisely, we have the following theorem which is a slight generalization of the result of  \cite{Chen}, see also \cite{Bol-Tre:ANTI}.
Take $E>\max_N W$. Suppose there exist a finite collection of nondegenerate  collision orbits $\{\gamma_k\}_{k\in K}$   of the billiard $(D,\partial D,H=E)$ starting and ending at the scatterer $N$, and define the graph $\Gamma$ as
in section \ref{sec:zero}.

\begin{thm}\label{thm:small}
There exists $\eps_0>0$ such that for any $\eps\in(0,\eps_0)$ and any path $\kbf=(k_j)_{j\in\Z}$ in the graph $\Gamma$
 there exists a unique trajectory of the billiard $(\Omega_\eps,\Sigma_\eps,H=E)$ shadowing  the collision chain $\gamma=(\gamma_{k_j})_{j\in\Z}$ of the degenerate billiard $(D,N\cup\partial D,H=E)$.
 This  trajectory is hyperbolic, and the set of all shadowing trajectories
form a hyperbolic invariant set $\Lambda_\eps$ for the billiard map $T_\eps:T^*\Sigma_\eps\to T^*\Sigma_\eps$.
\end{thm}

Thus we have a conjugacy  $\psi_\eps:\Sigma_\Gamma\to  \Lambda_\eps$ between $\calT:\Sigma_\Gamma\to\Sigma_\Gamma$ and $T_\eps:\Lambda_\eps\to\Lambda_\eps$.
The Lyapunov exponents of the trajectories in $\Lambda_\eps$ are large of order $|\ln\eps|$.
The shadowing error is of order $O(\eps)$, i.e.\  the shadowing orbit stays in a
$C\eps$-neighborhood of $\gamma$.
If the graph $\Gamma$ is branched (has two closed paths through the same vertex), then the billiard system  has a chaotic hyperbolic invariant set.

The proof of Theorem \ref{thm:small} uses the method of anti-integrable limit \cite{Aubry,Aubry-MacKay,Bol-Tre:ANTI}. In the next section we prove a more general result.

For the classical Birkhoff billiard in  a domain $D\subset \R^n$, 
Theorem \ref{thm:small} is proved in \cite{Chen}.
See also \cite{Bol-Tre:ANTI}.
If the boundary $\partial D$ is concave, then we have a dispersing Sinai billiard in $\Omega_\eps$,
so all trajectories are hyperbolic. In this case much stronger  results can be proved,
including ergodicity
\cite{Sinai:UMN,Sinai}.

\subsection{Billiards with thin scatterers}

\label{sec:thin}

We start with a Hamiltonian system $(M,H=E)$.
 Let $N$ be a submanifold in $M$   and $N_\eps $
its tubular $\eps$-neighborhood with the boundary $\Sigma_\eps=\partial N_\eps$.
Let $\Omega_\eps=M\setminus N_\eps$. Consider the billiard
$(\Omega_\eps,\Sigma_\eps,H=E)$ in $\Omega_\eps$
with the boundary $\partial \Omega_\eps=\Sigma_\eps$. As $\eps\to 0$, it  approaches the degenerate billiard $(M, N,H=E)$
with the scatterer $N$.

It is convenient to assume that $\Sigma_\eps$ is
defined as follows. Fix a Riemannian metric on $M$, not necessarily related to the Hamiltonian. For $x\in N$ let
$$
T_x^\perp N=\{u\in T_xM:u\perp T_xN\}.
$$
Let $f:T^\perp N\to M$, $q=f(x,u)$, be  the exponential  map $f(x,u)=\exp_x u$. Then  $f(x,0)=x$ and the  derivative  $D_u f(x,0):T_x^\perp N\to T_x^\perp N$ is an identity.
For any $x\in N$ take a convex neighborhood of $0$ in $T_x^\perp N$ containing 0
and let $S_x$ be its boundary. We assume that the hypersurface $S_x\subset T_x^\perp N$ is smooth, strictly convex (with positive definite second fundamental form)
 and smoothly depends on $x$, so that  $S=\cup_{x\in N}S_x\subset T^\perp N$ is a smooth manifold. Let
$$
\Sigma_\eps=f(\eps S)=\{q=f(x,\eps s):x\in N,\; s\in S_x\}.
$$
If $N$ is compact, or we take a compact subdomain in $N$, then for small $\eps>0$ $\Sigma_\eps$
is a smooth manifold and  $x\in N$, $s\in S_x$ are coordinates in $\Sigma_\eps$.

\begin{thm}\label{thm:thin-per}
Take $E>\max_NW$.
Let   $\gamma$ be a   nondegenerate periodic collision chain of the   the degenerate billiard $(M,N,H_0=E)$.
There exists $\eps_0>0$ such that for any $\eps\in (0,\eps_0)$ the   collision chain $\gamma$
is $O(\eps)$-shadowed by  a periodic orbit of the billiard in $\Omega_\eps$.  
\end{thm}

Recall that we call  the periodic collision chain $\gamma$ corresponding to $(\kbf,\xbf)$  nondegenerate if $\xbf$ is a nondegenerate  critical point of the  discrete action functional (\ref{eq:per}).
The shadowing periodic orbit has $\codim N$ large Lyapunov exponents of order $|\ln\eps|$.

A similar result holds for nonperiodic collision chains joining given points.

\begin{thm}\label{thm:thin-finite}
Let $\gamma$ be a finite nondegenerate collision chain joining the points $a,b\in M$.
There exists $\eps_0>0$ such that for any $\eps\in (0,\eps_0)$, $\gamma$  is
shadowed by  an orbit of the billiard in $\Omega_\eps$ joining the points $a,b$.
\end{thm}

Here nondegeneracy is understood in terms of the functional (\ref{eq:finite}).
Note that in Theorems \ref{thm:thin-per}--\ref{thm:thin-finite}, $\eps_0$ depends on the collision chain.
Under the hyperbolicity assumption  we can prove the following stronger theorem, where $\eps_0$ is independent of the trajectory.

Suppose that the degenerate  billiard satisfies appropriate conditions
so that the  corresponding DLS with the Lagrangian $\calL=\{L_k\}_{k\in K}$ has a compact hyperbolic $\calT$-invariant set $\Lambda\subset K^Z\times N^\Z$
of admissible trajectories $(\kbf,\xbf)$. Recall that a trajectory is admissible if the corresponding collision chain $\gamma=(\gamma_j)_{j\in\Z}$
satisfies the  jump condition $\Delta p_j\ne 0$ for all $j$.

\begin{thm}\label{thm:thin}
Suppose that the DLS   has  a compact hyperbolic invariant set $\Lambda$ of admissible trajectories.
There exists $\eps_0>0$ such that for $\eps\in (0,\eps_0)$ and any orbit $(\kbf,\xbf)\in\Lambda$
there exists a unique trajectory  of the billiard   $(\Omega_\eps,\partial \Omega_\eps,H=E)$
 which is $O(\eps)$-shadowing the corresponding collision chain $\gamma=(\gamma_{j})_{j\in\Z}$ of the degenerate billiard.
 This trajectory  is hyperbolic and
  all such trajectories form a compact hyperbolic invariant set $\Lambda_\eps\subset T^*\Sigma_\eps$.
 The billiard map  $T_\eps:\Lambda_\eps\to\Lambda_\eps$ is conjugate to the translation $\calT:\Lambda\to\Lambda$.
\end{thm}

Trajectories in $\Lambda_\eps$ have $\codim N$ large Lyapunov exponents of order $|\ln\eps|$.
A similar result holds for time periodic billiards.

\begin{rem} A more physical situation is when, as in section \ref{sec:small},  we start with an ordinary billiard $(D,\partial D,H)$
and delete a neighborhood $N_\eps$ of a submanifold $N\subset D$.
Let $\Omega_\eps=D\setminus N_\eps$. In the limit $\eps\to 0$  the billiard 
in the domain $\Omega_\eps$
will approach the degenerate billiard $(D,N\cup\partial D,H)$. The
 scatterer $\partial D\cup N$ has  components of different dimension:
the hypersurface $\partial D$ and a lower dimensional manifold $N$.  Theorem \ref{thm:thin} still holds.
However, for simplicity of notation we prove Theorem \ref{thm:thin} only for $D=M$ without boundary.
\end{rem}

A standard example   for Theorem \ref{thm:thin} is a system of $n$ small balls of radius $\eps$  moving
in a domain  $U\subset \R^d$
which  are elastically reflected  when colliding with the boundary or with each other.
Let  $U_\eps$ be $U$ with $\eps$-neighborhood of the
boundary deleted. Then
\begin{eqnarray*}
&\Omega_\eps=\{q=(q_1,\dots ,q_n)\in U_\eps^n: |q_i-q_j|\ge 2\eps \;\mbox{for $i\ne j$}\},\\
&\Sigma_\eps=\{q=(q_1,\dots ,q_n)\in U_\eps^n: |q_i-q_j|=\eps \;\mbox{for  some $i\ne j$}\},\\
&N= \{q=(q_1,\dots ,q_n)\in U^n: q_i=q_j \;\mbox{for  some $i\ne j$}\}.
\end{eqnarray*}

Under certain conditions on $U$ (for example, $U$ is a rectangular box) such billiard will be semi-dispersing  and then hyperbolicity and
ergodicity properties may be proved \cite{Sinai} for all trajectories of the billiard map.
Such results are important for the verification of the Boltzmann hypothesis \cite{Sinai,Simany}.  
Theorem \ref{thm:thin} implies much weaker results,
but it needs weaker assumptions. For example,  it does not need compactness
of the configuration space, so it can be used in situations studied in \cite{Kozlov}.

Next we give two simple concrete  examples.

\subsection{Examples}

\label{sec:discs}

Suppose two small balls of radius $\eps$ and masses $m_1,m_2$ move in the torus $\T^d$.
Then we have a billiard $(\Omega_\eps,\partial \Omega_\eps,H)$, where
\begin{equation}
\label{eq:torus}
\Omega_\eps=\{(q_1,q_2)\in\T^{2d}: |q_1-q_2|\ge 2\eps\},\qquad
H(q,p)=\frac{|p_1|^2}{2m_1}+\frac{|p_2|^2}{2m_2}.
\end{equation}
This is a billiard with thin scatterer, and the corresponding degenerate billiard is
$(\T^{2d},\Delta,H)$, where 
$$
\Delta=\{(q_1,q_2)\in\T^{2d}:q_1=q_2\}.
$$ 
This billiards has translational symmetry
and the corresponding momentum integral $p_1+p_2$, 
so the results of section \ref{sec:thin}
do not apply. However, we can reduce symmetry and obtain a reduced billiard $(\T^d\setminus B_\eps, \partial B_\eps, \tilde H)$, where $B_\eps$ is an $\eps$-ball and $\tilde H=|p|^2/2$.

The billiard $(\T^d\setminus B_\eps, \partial B_\eps, \tilde H)$ 
satisfies the conditions of Theorem \ref{thm:small}. The corresponding degenerate billiard  $(\T^d,\{0\}, \tilde H)$ is 0-dimensional. Fix energy $E>0$. The collision orbits $\gamma_k$ with end points at $0\in\T^d$ are labelled by the rotation vector $k\in\pi_1(\T^d)=\Z^d$. The graph $\Gamma$ has vertices in $\Z^d$ and two vertices $k,k'\in\Z^d$ are joined by an edge if $k\not\parallel k'$.
By Theorem \ref{thm:small},  for small $\eps>0$, any path in the graph is shadowed by a hyperbolic trajectory
of the billiard $(\T^d\setminus B_\eps, S_\eps, \tilde H=E)$.
This is true for any metric on a torus, only nongeneracy of geodesic loops is required. 

Of course the results provided by the  theory of dispersing billiards are much stronger.
Hyperbolicity and ergodicity of this  billiard was proved by Sinai \cite{Sinai:UMN}
for $d=2$, and then later for any $d\ge 2$ (see \cite{Sinai}).

As an example for Theorem \ref{thm:thin} consider  two balls of mass $m_1,m_2$ and radius $\eps$ freely moving in a  rectangle $U= (0,1)^d\subset\R^d$ with elastic reflections from the boundary and with each other. Let $U_\eps=(\eps,1-\eps)^d$. We obtain a billiard  in the domain
$$
\Omega_\eps=\{(q_1,q_2)\in U_\eps^2: |q_1-q_2|\ge2\eps\}.
$$
The Hamiltonian is as in  (\ref{eq:torus}), but now there is no momentum integral, so the Routh reduction is not possible. 
The corresponding degenerate billiard is $(U^2,\Delta,H)$.
For fixed energy $E=1/2$, the Jacobi metric in $U^2$ is
 $$
\|\dot q\|=\sqrt{m_1|\dot q_1|^2+m_2|\dot q_2|^2}.
$$

 The Lagrangians  
 of the corresponding DLS are functions on $U^2$ defined as follows. 
 Take two points $q_\pm\in U$ and two billiard trajectories $\gamma_1,\gamma_2$ of an ordinary billiard  in $U$ joining these  points after one or several reflections from the boundary $\partial U$.   Let $l_1,l_2$ be  their lengths.
 There are an infinite number of such pairs of billiard trajectories, we label them by an index $k$.
 Then the $k$-th Lagrangian is
 $$
 L_k(q_-,q_+)=\int \sqrt{m_1|\dot \gamma_1|^2+m_2|\dot \gamma_2|^2}\,dt
 = \sqrt{m_1l_1^2+m_2l_2^2},\qquad q_\pm \in U.
$$
The Lagrangians $L_k$ are (nonstrictly) convex on $U^2$.

For collision chains with fixed end points $a,b\in U$, the  action functional (\ref{eq:finite})
is strictly convex. Thus any such collision chain is nondegenerate, and by Theorem \ref{thm:thin-finite} for small $\eps>0$ it is shadowed by a trajectory of the billiard in $\Omega_\eps$.

One can check that also the periodic action functional (\ref{eq:per}) is strictly convex
for certain types of  periodic codes $\kbf$.  Then we can apply Theorems \ref{thm:thin-per} and  obtain shadowing periodic orbits. 

The simplest situation is for $d=1$, when collision chains correspond
to billiard trajectories in a triangle (see e.g.\ \cite{Koz-Tre:Bill}). Then periodic trajectories
with odd number of reflections are nondegenerate (degenerate for even number of reflections).

Theorem \ref{thm:thin} does not apply in this example because the  degenerate billiard has no hyperbolic invariant sets.  Thus our approach does not provide chaotic hyperbolic invariant sets for the billiard in $\Omega_\eps$.
The methods of the theory of semi-dispersing billiards \cite{Sinai}
make it possible to prove much stronger results, including hyperbolicity and ergodicity, see e.g.\  \cite{Simany}. 

Applications of our method makes more sense for nondispersing billiards.
For example, two small balls moving inside a spherical domain.

\subsection{Proof of Theorem \ref{thm:thin}}

\label{sec:proof}

We will prove only Theorem \ref{thm:thin}, the proofs of Theorems \ref{thm:thin-per}--\ref{thm:thin-finite} are similar but simpler.

We represent any point $q=f(x,\eps s)$ on the boundary $\Sigma_\eps$ by a pair $x\in N$, $s\in S_x\subset T_x^\perp N$.
Then the configuration space of the billiard $(\Omega_\eps,\Sigma_\eps,H=E)$ is identified with
the manifold
$$
Q=\{q=(x,s): x\in N,\; s\in S_x\},
$$
and the billiard map $T_\eps$ is a map of a subset in $T^*Q$.

Let $\calL=\{L_k\}_{k\in K}$ be the discrete Lagrangian for the degenerate billiard $(M,N,H=E)$.
Then $L_k$ is a function on an open set $U_k\subset N^2$ defined as the Maupertuis action $L_k(x_-,x_+)=J_E(\gamma)$ of a trajectory
$\gamma=\gamma(x_-,x_+)$ with energy $E$ joining $x_-$ and $x_+$.

The  discrete Lagrangian $L_{k}^\eps(q_-,q_{+})$, $q_\pm=(x_\pm,s_\pm)$, for the billiard $(\Omega_\eps,\Sigma_\eps,H=E)$ is
 the local generating function of the billiard map $T_\eps$. This is a function on an open subset
 $$
 U_k^\eps=\{(q_-,q_+)\in Q\times Q: q_\pm=(x_\pm,s_\pm),\; (x_-,x_+)\in U_k,\; s_\pm \in S_{x_\pm}\}
 $$
 defined as the Maupertuis action $J_E(\gamma^\eps)$ of the trajectory $\gamma^\eps=\gamma^\eps(x_-,s_-,x_+,s_+)$ of energy $E$
 joining the points $q_-(\eps)=f(x_-,\eps s_-)$ and $q_+(\eps)=f(x_{+},\eps s_{+})$.
 Since $x_\pm$ are non-conjugate, $\gamma^\eps$
 smoothly depends on $\eps$.  We have
$$
\frac{\partial}{\partial\eps}\bigg|_{\eps=0}q_\pm(\eps)=s_\pm.
$$
Hence the variation of the Maupertuis action is
$$
\frac{\partial}{\partial\eps}\bigg|_{\eps=0}J_E(\gamma^\eps)=\frac{\partial}{\partial\eps}\bigg|_{\eps=0}\int_{\gamma^\eps}p\, dq=\langle p_+,s_+\rangle - \langle p_-,s_-\rangle.
$$
where $p_\pm(x_-,x_+)$ are the initial and final momenta of $\gamma$. We obtain the generating function of the billiard map $T_\eps$ in the form
\begin{eqnarray*}
L_{k}^\eps(q_-,q_+)&=&J_E(\gamma^\eps)\\
&=&L_{k}(x_-,x_{+})-\eps \langle p_-(x_-,x_+),s_-\rangle +\eps \langle p_+(x_{-},x_{+}),s_{+}\rangle+O(\eps^2).
\end{eqnarray*}

 The formal discrete action functional for the billiard $(\Omega_\eps,\Sigma_\eps,H=E)$  has the form
\begin{equation}
\label{eq:Aeps}
A_\kbf^\eps(\qbf)=\sum_{j\in\Z} L_{k_j}^\eps(q_j,q_{j+1})=A_\kbf(\xbf)+\eps \sum_{j\in\Z} \langle \Delta p_j,s_j\rangle +O(\eps^2),
\end{equation}
where
\begin{eqnarray*}
\qbf=(\xbf,\sbf)\in Q^\Z,\qquad (x_j,x_{j+1})\in U_{k_j},\quad s_j\in S_{x_j},\\
\Delta p_j(\kbf,\xbf)=p_j^--p_j^+,\qquad p_j^-=p_-(x_j,x_{j+1}),\quad p_j^+=p_+(x_{j-1},x_{j}).
\end{eqnarray*}
As above, the functional $A_\kbf^\eps$ is formal but its derivative $DA_\kbf^\eps(\qbf)$ is well defined, so $O(\eps^2)$ means such a term in the derivative.

For fixed $(\kbf,\xbf)$, the functional (\ref{eq:Aeps}) splits in the sum of independent functions
$$
g_j(\kbf,\xbf,s_j)=\langle \Delta p_j(\kbf,\xbf),s_j\rangle,\qquad s_j\in S_{x_j}.
$$
Since the surface $S_{x_j}$ is strictly convex with curvature bounded away from 0,
the function $s_j\to g_j(\kbf,\xbf,s_j)$ has a unique nondegenerate maximum point $s_j(\kbf,\xbf)$ such that $s_j\perp \Delta p_j$. Let $\sbf(\xbf)=(s_j(\kbf,\xbf))_{j\in\Z}$.

Suppose that $\Lambda\subset K^\Z\times N^\Z$ is a compact hyperbolic invariant set of admissible trajectories of the DLS describing the degenerate billiard.
For any $(\kbf,\xbf)\in \Lambda$ and $j\in\Z$ let $\gamma_j:[t_j,t_{j+1}]\to M$
be the corresponding  collision orbit with energy $E$ joining the points $x_j,x_{j+1}\in N$.
Then the jump of momentum at $j$-th collision is $\Delta p_j=\Delta p_j(\kbf,\xbf)$.

Since $\Lambda$ is compact,  the jump condition  is uniform :
$\|\Delta p_j(\kbf,\xbf)|\ge \delta>0$  with $\delta=\delta(\Lambda)>0$ independent of the trajectory  in $\Lambda$. Hence
 \begin{equation}
 \label{eq:convex}
 \|(D^2_{s_j}g_j(\kbf,\xbf,s_j))^{-1}\|\le C_0=C_0(\Lambda).
 \end{equation}

By the definition of a hyperbolic invariant set, $\xbf$ is a nondegenerate critical point of $A_\kbf$, and moreover $A''_\kbf(\xbf)$
 has bounded inverse:
 \begin{equation}
 \label{eq:Fk}
 \|(A''_\kbf(\xbf))^{-1}\|_\infty\le C_1=C_1(\Lambda)
 \end{equation}
 with $C_1$ independent of the trajectory. We have
\begin{eqnarray*}
&D_\xbf A_\kbf^\eps(\qbf)=F_\kbf(\xbf)+\eps R(\kbf,\qbf,\eps),\qquad &F_\kbf(\xbf)=A'_\kbf(\xbf),\\
&D_\sbf A_\kbf^\eps(\qbf)=\eps H_\kbf(\qbf)+\eps^2 S(\kbf,\qbf,\eps),\qquad &H_\kbf(\qbf)=(D_{s_j}g_j(\kbf,\xbf,s_j))_{j\in\Z},
\end{eqnarray*}
where
$$
\|D_\qbf R(\kbf,\qbf,\eps)\|\le C_2,\quad \|D_\qbf S(\kbf,\qbf,\eps)\|\le C_2,\quad
\|D_\qbf H_\kbf(\qbf)\|\le C_2,
$$
where $C_2=C_2(\Lambda)$.
By (\ref{eq:convex}),
\begin{equation}
\label{eq:Hk}
\|(D_\sbf H_\kbf(\qbf))^{-1}\|_\infty\le C_0.
\end{equation}

 By (\ref{eq:Fk})--(\ref{eq:Hk}),  there exists $\eps_0>0$
 such that for $\eps\in (0,\eps_0)$ the functional $A_\kbf^\eps$ has a nondegenerate critical
point $\qbf$ close to $(\xbf,\sbf(\xbf))$ such that
$$
\|(D^2A_\kbf^\eps(\qbf))^{-1}\|_\infty\le \eps^{-1}C(\Lambda).
$$
Hence $(\kbf,\qbf)$  corresponds to a hyperbolic trajectory of the billiard  $(\Omega_\eps,\Sigma_\eps,H=E)$.
Theorem \ref{thm:thin} is proved.
\qed

\medskip

When $N$ is a discrete set, there is no variable $\xbf$ and we obtain the proof of Theorem \ref{thm:small}.
In this case the functional $A_\kbf^\eps(\sbf)$ has an anti-integrable form \cite{Aubry}, see \cite{Chen,Bol-Tre:ANTI}.

\section{Systems with Newtonian singularities}
\label{sec:sing}

In this  section we discuss a similar but technically more difficult example of  degenerate billiards
arising in applications to Celestial Mechanics. The proofs are postponed to the next publication \cite{Bol:Bill_RCD}.

Consider a Hamiltonian system $(M\setminus N,H_\mu)$ on $T^*(M\setminus N)$ with a classical smooth\footnote{$C^4$ is enough.}
Hamiltonian
\begin{equation}
\label{eq:Hmu}
H_\mu(q,p)=\frac12\|p-w_\mu(q)\|_\mu^2+W_\mu(q)+\mu V(q,\mu)
\end{equation}
depending on a small parameter $\mu\in (-\mu_0,\mu_0)$. Here $\|\;\|_\mu$ is a Riemannian metric on $M$ smoothly depending on the parameter $\mu$, and $w_\mu$ and $W_\mu$ are covector field and a function on $M$ smoothly depending on $\mu$. The potential $V$ is smooth on $M\setminus N$ but undefined on $N$.

For $\mu=0$ the unperturbed system $(M,H_0)$ is a classical  Hamiltonian system on $T^*M$
with Hamiltonian $H_0$ of the form (\ref{eq:HL}). The perturbation consists of two parts: regular perturbation
which is a smooth function on $T^*M$, and a singular part $\mu V$.

We say that $V$ has a Newtonian singularity on $N$
if in a tubular neighborhood of $N$
there exists a smooth nonzero function $\phi$  such that
\begin{equation}
\label{eq:V}
V(q,\mu)=-\frac{\phi(q,\mu)}{d(q,N)}.
\end{equation}
The distance $d$ is defined by the Riemannian metric $\|\;\|_\mu$. For definiteness we assume $\mu>0$. Then if $\phi<0$, the singular force $-\mu\nabla V$ is repelling (like the Coulomb force),
and if $\phi>0$ attracting (like the gravitational force).

 We are interested in nearly collision trajectories $\gamma_\mu:\R\to M\setminus N$ of system $(M\setminus N,H_\mu)$ which come $O(\mu)$-close to $N$.
 Their  limits  as $\mu\to 0$ are collision chains of the degenerate billiard $(M,N,H_0)$ with Hamiltonian $H_0$ and scatterer $N$.

 \begin{rem}
 One can consider other  (nonphysical) singularities replacing $d(q,N)$
 in (\ref{eq:V}) with $d^\alpha(q,N)$, $\alpha>0$.
 However, then our methods  do not work.
 \end{rem}

 \subsection{General $n$ center problem}

 \label{sec:cent}

 For simplicity first  let $N=\{a_i\}_{i\in J}$ be a finite set.
 When $\phi>0$ (attracting force), we
call the system $(M\setminus N,H_\mu)$ the $n$-center
problem.

 We fix energy $H_\mu=E>\max_NW_0$  and study the system $(M\setminus N,H_\mu=E)$  on the energy level.
Let $(M,N,H_0=E)$ be the  0-dimensional degenerate billiard with the scatterer $N$.
Suppose there are several nondegenerate collision orbits
$\{\gamma_k\}_{k\in K}$ connecting points
$a_k^-,a_k^+\in N$.  Let $p_{k}^-$ and
$p_k^+$ be the initial and final momenta of $\gamma_k$.
Define an oriented  graph $\Gamma$ as in section \ref{sec:zero}   joining $k$,
$k'$ by an edge if $a_k^+=a_{k'}^-$ and $p_{k}^+\ne p_{k'}^-$.

\begin{thm}
\label{thm:centers}
There exists $\mu_0>0$ such that for all
$\mu\in(0,\mu_0)$ and any path $\kbf=(k_i)_{i\in\Z}$
in the graph $\Gamma$ there exists a unique (up to a time shift)
trajectory of system $(M\setminus N,H_\mu=E)$ which is $O(\mu)$-shadowing the collision chain
$\gamma=(\gamma_{k_i})_{i\in\Z}$ of the degenerate billiard.
This trajectory is hyperbolic.
\end{thm}

Hence there is a hyperbolic  invariant subset in $\{H_\mu=E\}\subset T^*(M\setminus N)$
on which the
system is  conjugate to a suspension of a  topological Markov chain.
The topological entropy is positive provided that the graph
$\Gamma$ has a connected branched subgraph.

\begin{rem}
In the attracting case $\phi>0$ the shadowing trajectory may have
regularizable collisions with the singular set $N$.
Thus dynamics is well defined. However, if we want to avoid collisions,
we should add  to the definition of an edge $(k,k')$ of the graph $\Gamma$ the condition  that $v_{k}^+\ne -v_{k'}^-$,
where $v_k^\pm$ are the initial and final velocities  of $\gamma_k$.
This condition is less essential than the changing direction condition $p_{k}^+\ne p_{k'}^-$.
Note that $v_{k}^+\ne -v_{k'}^-$ is not equivalent to $p_{k}^+\ne -p_{k'}^-$ if the system
is nonreversible.
\end{rem}

The proof of Theorem \ref{thm:centers} is given in \cite{Bol-Mac:centers} (using different language) for $d=\dim M=2$ by
using the Levi-Civita regularization of collisions and the Shilnikov Lemma \cite{Shiln}.
Another ingredient is the method of an anti-integrable limit \cite{Aubry}.
In \cite{Bol-Mac:spatial} the proof was given
for $d=3$ by using the KS regularization.
In these papers only the gravitational case was considered, but
in the repelling case the proof is the same.

\subsection{Classical $n$ center problem}

\label{sec:cent-class}

The classical example is when $M=\R^d$ with the Euclidean metric and
$$
H_\mu(q,p)=\frac12 |p|^2+\mu V(q),\qquad V(q)=-\sum_{i=1}^n\frac{\alpha_i}{|q-a_i|},\qquad \alpha_i>0.
$$
Note that  small $\mu>0$ is equivalent to large energy $E>0$ (the time change $t\to t\sqrt{\mu}$
replaces $\mu\to 1$ and $E\to E/\mu$).
Trajectories of the corresponding degenerate billiard are polygons with vertices $a_i$.
By Theorem \ref{thm:centers}, for small $\mu>0$ (or large energy) any path in the polygon is shadowed
by a trajectory of the $n$-center problem.
This case was studied in \cite{Knauf}. In fact for $d=2$ the curvature of the Jacobi metric is negative, so stronger results can be obtained.
This does not work when $d\ge 3$, when some sectional curvatures are negative.

Theorem \ref{thm:centers} provides chaotic trajectories only for $n>3$,
and when not all centers are collinear (then the graph $\Gamma$ is branched).
In fact the $n$ center problem has chaotic trajectories for any $n\ge 3$ for purely topological reasons. Moreover smallness of $\mu>0$ (or large energy) is not needed, and one can add to the potential a smooth negative function. The proof is given
in \cite{Bol:centers} for $d=2$ and  in \cite{Bol-Neg:reg} for $d=3$.
For $d=3$ the proof is nontrivial: it is based on the global KS regularization and deep
results of Gromov and Paternain.
Probably this result is true also for $d\ge 4$, but a different proof is needed.
The approach based on degenerate billiards is much simpler, but it it restricted to the
case of small $\mu>0$.

\subsection{Restricted 3 body problem}

\label{sec:restrict}

A more important classical example is the plane restricted
circular  3 body problem: a small body (Asteroid) moves in $\R^2$ under the
action of gravitational forces of Jupiter with mass $\mu$  and the Sun with mass $1-\mu$. Suppose the  Sun is at the origin and
Jupiter is moving along a unit circle.  Let $U=\R^2\setminus\{0\}$.
In the rotating coordinate frame we obtain a Hamiltonian $H_\mu=H_0+\mu V$ of the form
(\ref{eq:Hmu}) with
 $$
H_0(q,p)=\frac12|p|^2+ q_1p_2-q_2p_1 -\frac{1}{|q|},\qquad
  V(q)=-\frac 1{|q-a|},\quad a=(1,0)\in\R^2.
$$
This is a Hamiltonian system $(U\setminus \{a\},H_\mu)$ with a Newtonian singularity at $a$.
In the limit $\mu\to 0$ we obtain a 0-dimensional degenerate billiard $(U,\{a\},H_0)$, where
$H_0$ is the Hamiltonian  of the Kepler problem in a rotating coordinate frame.
This is an integrable system, but a description of collision orbits with given Jacobi integral $H_0=E$
is not easy since it involves
solving the transcendental Kepler's equation \cite{AKN}.

However one can prove \cite{Bol-Mac:centers} that for values of $E$ such that the Kepler ellipse can cross the unit circle (Jupiter's orbit),
the graph $\Gamma$ describing collision chains has  an infinite number of vertices
corresponding to the number of revolutions $k\in \Z$ between collisions
By Theorem \ref{thm:centers}, hyperbolic dynamics of nearly collision trajectories is quite rich. 
Such trajectories are called second species solutions of Poincar\'e. The existence of two link periodic second species solutions was proved much earlier in \cite{Mar-Nid}.

Second species solutions have  a long history in Astronomy (see the bibliography in \cite{Henon}), 
but without formal proofs.
For example, the trajectory of Voyager is a second species solution.

If Jupiter's orbit is not a circle but an ellipse, we obtain the elliptic restricted 3 body problem.
Then the corresponding degenerate billiard is time periodic, and the scattering map is a multivalued twist map of a cylinder.
Then there is fast chaotic  ``diffusion'' of the Jacobi constant.
For details see \cite{Bol:Nonlin}.

\subsection{General systems with Newtonian singularities}

\label{sec:newton}

Now consider a general Hamiltonian system (\ref{eq:Hmu}) with $N$ a submanifold  in $M$ and
the corresponding DLS
describing the degenerate billiard $(M,N,H_0=E)$.  

\begin{thm}\label{thm:per}
Let   $\gamma$ be a nondegenerate periodic  chain of the   of the degenerate billiard $(M,N,H_0=E)$.
There exists $\mu_0>0$ such that for any $\mu\in (0,\mu_0)$ the  collision chain $\gamma$
is shadowed by  a periodic orbit of the system $(M\setminus N,H_\mu=E)$.
\end{thm}

The shadowing error is of order $O(\mu^\alpha)$ for any $\alpha\in (0,1)$.
Under small additional assumptions, this can be improved to $O(\mu|\ln \mu|)$.
Theorem \ref{thm:per} is a generalization of a theorem in \cite{Bol-Neg:RCD}, 
where it was proved for the plane 3 problem,
see  section \ref{sec:cel}.

\begin{thm}
\label{thm:hyp}
Let $\Lambda\subset K^\Z\times N^\Z$ be a compact hyperbolic invariant set of the DLS   such that all orbits in $\Lambda$ are admissible.
There exists $\mu_0>0$ such that for any $\mu\in (0,\mu_0)$
and  any orbit $(\kbf,\xbf)\in\Lambda$
there exists an trajectory of system $(M\setminus N,H_\mu=E)$  shadowing
the corresponding collision chain of the degenerate billiard. Such trajectories form a compact hyperbolic invariant set $\Lambda_\mu\subset \{H_\mu=E\}\subset T^*(M\setminus N)$ of system
$(M\setminus N,H_\mu=E)$.
\end{thm}

The shadowing error is of the same order  as in Theorem \ref{thm:per}.
Note that in  Theorem \ref{thm:per} the periodic orbit  
does not need to be hyperbolic, so Theorems \ref{thm:per} and \ref{thm:hyp} are formally independent.
We will prove both theorems in the  next paper \cite{Bol:Bill_RCD}.
Also an analog of Theorem \ref{thm:thin-finite} holds.

\begin{rem}
Formally these results do not work for systems with symmetry in Celestial Mechanics,
see   section \ref{sec:cel}.
Indeed, systems in Celestial Mechanics have continuous symmetries,
so  any critical point of the action functional is degenerate.
One can prove an analog of Theorem \ref{thm:per} for systems with symmetry (this is proved in
\cite{Bol-Neg:DCDS} for a particular case).
Another option is to perform  the Routh reduction eliminating symmetry.
Then also Theorem \ref{thm:hyp} can be applied.
\end{rem}

\begin{rem}
As in the $n$-center problem,
for $\phi>0$ (attracting  singularity) we   need an extra condition to avoid regularizable singularities of the shadowing orbit.
Let $v_j^\pm=H_p(x_j,p_j^\pm)$ be the velocity at $j$-th collision and let $u_j^\pm$ be its projection to the quotient space $T_{x_j}M/T_{x_j}N$.
We assume $u_j^+\ne -u_j^-$ (no straight reflection). As noted in the previous section, this condition is less essential
since  singularities are regularizable, so dynamics is always well defined,
also for trajectories colliding with $N$.
\end{rem}

\subsection{Second species solutions of the 3 body problem}

\label{sec:cel}

As an example, consider the plane 3-body problem with masses $m_1,m_2,m_3$, where
 $m_3$ is much larger than $m_1,m_2$:
$$
\frac{m_1}{m_3}=\mu\alpha_1,\quad \frac{m_2}{m_3}=\mu\alpha_2,\quad \alpha_1+\alpha_2=1,\quad \mu\ll 1.
$$

Let $p_1,p_2,p_3\in \R^2$ be the momenta of the bodies. Assume that the center of
mass is at rest: $p_1+p_2+p_3=0$. Let $q_1,q_2\in\R^2$ be positions of $m_1,m_2$ relative to $m_3$.
Set $q=(q_1,q_2)$, $p=(p_1,p_2)$.  After rescaling, the motion of $m_1,m_2$ relative to $m_3$ is
described by the Hamiltonian
\begin{equation}
 H_\mu(q,p)=H_0(q,p)+\mu \frac{|p_1+p_2|^2}2 -
\mu\frac{\alpha_1\alpha_2}{|q_1-q_2|},
\end{equation}
where
\begin{equation}
\label{eq:H0}
H_0=H_1+H_2,\qquad H_i= \frac{|p_i|^2}{2 \alpha_i}-\frac{\alpha_i}{|q_i| }.
\end{equation}
The Hamiltonian $H_\mu$ has the form (\ref{eq:Hmu}) with
$$
V(q)=-\frac{\alpha_1\alpha_2}{|q_1-q_2|}.
$$
The configuration  space  is $U^2\setminus\Delta$, where $U=\R^2\setminus\{0\}$ and
$$
\Delta=\{q=(q_1,q_2)\in U^2: q_1=q_2\}.
$$

In the limit $\mu\to 0$ the small bodies   $m_1$ and $m_2$
do not interact and the Hamiltonian system  $(U^2\setminus \Delta,H_\mu)$ becomes a degenerate billiard $(U^2,\Delta,H_0)$
with the scatterer $\Delta$.
The unperturbed Hamiltonian $H_0$ describes 2
uncoupled Kepler problems with  masses $\alpha_i$. It has integrals of energy $H_i$
and   of angular momenta $G_i$
corresponding to the rotational symmetry.\footnote{Since the Kepler problem is totally degenerate,
there are  also Laplace integrals.}
For $\mu>0$ the 3 body problem has integrals of energy  $H_\mu=E$ and total angular momentum
$G=G_1+G_2$, but no more.\footnote{Formally this is not proved for all values of parameters,
but  there is no doubt that this is true in general.}

Collision chains of the degenerate billiard $(U^2,\Delta,H_0)$ are concatenations  of pairs of arcs of Kepler orbits.
By the reflection law (\ref{eq:Deltap})--(\ref{eq:DeltaH}), the total energy and total momentum
\begin{eqnarray*}
\label{eq:energy}
H_0&=&\alpha_1h_1+\alpha_2h_2,\qquad h_i=\frac12|\dot q_i|^2-\frac{1}{|q_i|},\\
 y&=&p_1+p_2=\alpha_1\dot q_1+\alpha_2\dot q_2
 \end{eqnarray*}
are continuous at collisions.
Then the total angular momentum $G=G_1+G_2$ and total energy $H_0=E$ are
constant along a collision chain, but not the individual energies
$H_1,H_2$ and angular momenta $G_1,G_2$ of the bodies $m_1,m_2$.

For  trajectories of the 3-body problem shadowing collision chains  the small bodies move nearly along Kepler ellipses
and after many revolutions they almost collide.
Near collision, they move around their common center of mass nearly along small Kepler hyperbolas.
After exiting a neighborhood of a collision, they move along a new pair of Kepler ellipses until new nearly collision and so on.

Trajectories of the 3-body problem shadowing collision chains were named by Poincar\'e 
the second species solutions.
Poincar\'e didn't provide a rigorous  proof of the existence of the second species solutions.
There were many works of Astronomers (see the bibliography in \cite{Henon}),
but mostly they studied (without formal proofs) solutions which come $O(\mu^\alpha)$, $0<\alpha<1$,
close to collisions, while shadowing trajectories we study come $O(\mu)$-close to collisions.
The existence of $O(\mu)$-shadowing trajectories was proved in \cite{Bol-Neg:RCD}
for the case $\alpha_1\ll\alpha_2$ by using a version of Theorem \ref{thm:per}. The general
case will be studied in a subsequent publication.

We briefly describe the  discrete Lagrangians $\calL=\{L_k\}_{k\in\Z^2}$
of the degenerate billiard $(U^2,\Delta,H_0=E)$.
A collision orbit $\gamma:[t_-,t_+]\to U^2$  is a pair
$(\gamma_1,\gamma_2)$ of Kepler orbits
$\gamma_i:[t_-,t_+]\to U$ with energies $h_i$ such that
$$
\gamma_1(t_\pm)=\gamma_2(t_\pm)=x_\pm,\qquad \alpha_1h_1+\alpha_2h_2=E.
$$
Let $k_i$ is the number of revolutions of $i$-th body between collisions and $k=(k_1,k_2)\in\Z^2$.
The Maupertuis action of $\gamma$ defines the discrete Lagrangian
 \begin{equation}
\label{eq:A}  L_k(x_-,x_+)=J_E(\gamma)= \int_\gamma p_1\, dq_1+p_2\,
dq_2=\alpha_1 J_{h_1}(\gamma_1)+\alpha_2 J_{h_2}(\gamma_2),
\end{equation}
where
\begin{equation}
\label{eq:Ih} J_{h_i}(\gamma_i)=\int_{\gamma_i}
\sqrt{2|x|^{-1}+2h_i}\,|dx|
\end{equation}
is the Maupertuis action of a Kepler orbit with energy $h_i$.
A more explicit formula   is given in
\cite{Bol-Neg:DCDS}:
$$
L_k(z)=\min_{\alpha_1h_1+\alpha_2h_2=E}(\alpha_1
J_{k_1}(h_1,z)+\alpha_2 J_{k_2}(h_2,z)),\qquad z=(x_-,x_+),
$$
where
$$
J_n(h,z)=(-2h)^{-1/2}(2\pi|n|+(\sgn n)f(-2hz)),
$$
and $f(x_-,x_+)$ is the action of a simple Kepler arc with major semiaxis 1 and mass 1
joining $x_\pm$.

We also need to avoid early collisions of $\gamma_1,\gamma_2$ for $t_-<t<t_+$. They happen if the periods are commensurable:
$$
n_1T_1=n_2T_2,\qquad 0<n_i<k_i,\qquad T_i=2\pi (-2h_i)^{-3/2}.
$$
For large $\|k\|$  almost all collision orbits do not satisfy these  equations.

Another difficulty is that the  discrete billiard
has the integral of angular momentum. Thus to apply Theorem \ref{thm:hyp}
we need to  reduce rotational symmetry by the Routh method (see  \cite{Bol-Neg:DCDS}).
The details  will be given in  a subsequent publications.

\end{document}